\documentclass{icmart}
\usepackage{amssymb}
\usepackage{amsmath}
\usepackage[pdftex]{color,graphicx}

\contact[jean-francois.legall@math.u-psud.fr]{Universit\'e Paris-Sud, Math\'ematiques, B\^at.425, 
91405 Orsay C\'edex, France}

\newtheorem{theorem}{Theorem}[section]
\newtheorem{lemma}[theorem]{Lemma}
\newtheorem{proposition}[theorem]{Proposition}
\newtheorem{corollary}[theorem]{Corollary}

\newtheorem{definition}[theorem]{Definition}

\def\build#1_#2^#3{\mathrel{ \mathop{\kern 0pt#1}\limits_{#2}^{#3}}}

\title[Random Geometry on the Sphere]{Random Geometry on the Sphere}

\author[Jean-Fran\c cois Le Gall]{Jean-Fran\c cois Le Gall\thanks{This paper will appear in the Proceedings 
of ICM Seoul 2014}}

\begin{document}

\begin{abstract}
We introduce and study a universal model of random geometry in two dimensions. 
To this end, we start from a discrete graph drawn on the sphere, which is chosen uniformly
at random in a certain class of graphs with a given size $n$, for instance the class
of all triangulations of the sphere with $n$ faces. We equip the vertex set of the graph
with the usual graph distance rescaled by the factor $n^{-1/4}$. We then prove that the resulting random metric space
converges in distribution as $n\to\infty$, in the Gromov-Hausdorff sense, toward a limiting
random compact metric space called the Brownian map, which is universal in the sense that it
does not depend on the class of graphs chosen initially.  The Brownian map is homeomorphic to the sphere,
but its Hausdorff dimension is equal to $4$. We obtain detailed information about the structure of
geodesics in the Brownian map. We also present the infinite-volume variant of the Brownian
map called the Brownian plane, which arises as the scaling limit of the uniform infinite planar quadrangulation.
Finally, we discuss certain open problems. This study is motivated in part by
the use of random geometry in the physical theory of two-dimensional quantum gravity.
\end{abstract}

\begin{classification}
05C80, 60D05
\end{classification}

\begin{keywords}
Planar map, Gromov-Hausdorff convergence, Brownian map.
\end{keywords}

\maketitle

\section{Introduction}

In the last ten years, there has been much interest in discrete and continuous models of random geometry
in two dimensions. As a first naive attempt to construct a continuous model, one could imagine choosing
at random a Riemannian metric on the sphere. However there seems to be no canonical way of making such a
choice, and a better approach leading to a continuous object that is universal in some sense is to start
from discrete models. Informally, we will consider a large graph drawn on the sphere and chosen at random in
a suitable class: We then expect that the suitably rescaled graph distance on the vertex set will converge in an appropriate sense, when the size of the graph tends to infinity,
to a random metric on the sphere.

 More precisely, we consider planar maps, which are finite connected graphs embedded
in the sphere $\mathbb{S}^2$ -- more precise definitions will be given in Section 2 below. One is interested in the 
``shape'' of the graph, so that two planar maps are identified if the second one is the image of
the first one under an orientation-preserving homeomorphism of the sphere. The faces of the 
map are the connected components of the complement of edges, and a planar map is called a
triangulation if all faces are triangles, possibly with two edges glued together. For technical reasons, it is
convenient to deal with {\it rooted} triangulations, meaning that there is a distinguished oriented edge called
the root edge, whose tail is the root vertex. Thanks to the
preceding  identification, for every fixed even integer $n$, there are only a finite number of rooted triangulations with 
$n$ faces, and therefore it makes sense to choose one of them at random, which we denote by $T_n$. 
Then, supposing that there is a canonical way of embedding $T_n$ in the sphere (recall that
the embedding of a planar map is only defined up to orientation-preserving homeomorphisms), it seems reasonable
to conjecture that the vertex set $V(T_n)$ of $T_n$ will become dense in $\mathbb{S} ^2$ when $n\to\infty$, and that the graph
distance $d_{\rm gr}$ on this vertex set will converge, modulo a suitable rescaling, to  a random metric on $\mathbb{S}^2$.

The preceding program is still the subject of active research -- there are indeed (almost) canonical ways of
embedding triangulations, using circle packings or the uniformization theorem for Riemann surfaces, see the discussion in
Section \ref{secopen} below.
Here we will adopt a slightly different point of view. We consider the finite metric space $(V(T_n),d_{\rm gr})$
as a (random) element of the set $\mathbb{K}$ of all isometry classes of compact metric spaces. The space
$\mathbb{K}$ is equipped with the Gromov-Hausdorff distance $d_{GH}$ (cf. subsection \ref{GHsubsec}), and $(\mathbb{K},d_{GH})$
is both separable and complete. A key result \cite{LGU}, which solves a conjecture of Oded Schramm \cite{Sch}, states that
\begin{equation}
\label{convtri}
(V(T_n),6^{1/4}n^{-1/4}d_{\rm gr})\build{\longrightarrow}_{n\to\infty}^{\rm(d)} ({\bf m}_\infty, D),
\end{equation}
where the convergence holds in distribution in the space $(\mathbb{K},d_{GH})$, and the limit $({\bf m}_\infty, D)$
is a random metric space, which is called the Brownian map after Marckert and Mokkadem \cite{MaMo}. Despite the
somewhat unusual fact that
we are dealing with random compact metric spaces, (\ref{convtri}) is nothing but a particular case of
the standard notion of convergence in distribution for random variables with values in a Polish space. Note that the constant
$6^{1/4}$ in (\ref{convtri}) is just a normalization factor. 

A very important feature of the convergence (\ref{convtri}) is the fact that it holds for
much more general random planar maps, with the {\it same} limiting space $({\bf m}_\infty, D)$ up to unimportant scaling constants. Recall that a planar map is 
a $p$-angulation if all faces have degree $p$. Then it was proved in \cite{LGU} that an analog of (\ref{convtri}) holds for $p$-angulations with a fixed number of faces, for every even integer $p\geq 4$.
Note that the 
special case of quadrangulations, corresponding to $p=4$, has been treated independently by Miermont \cite{Mi-quad}
via a different method. Similarly, analogs of \eqref{convtri} hold for general planar maps with a fixed number of edges \cite{BJM} or
for bipartite planar maps with a fixed number of edges \cite{Ab}. It is also possible to impose local constraints on the planar maps: 
A result similar to (\ref{convtri}) holds for simple triangulations or quadrangulations \cite{AA}, where no loops or multiple edges are allowed, or for
quadrangulations with no pendant vertices \cite{BLG}. One indeed expects that the Brownian map will appear as
the scaling limit of very general random planar maps provided some bound holds on the distribution of the degree of 
a typical face (on the other hand, the paper \cite{LGM2} shows that different scaling limits may occur
if one considers distributions that favour the appearance of ``very large'' faces).

The Brownian map thus appears as a {\it universal} object, in the sense that it is the scaling limit of many different
models of random planar maps: This is of course similar to the case of Brownian motion, which is the universal scaling 
limit of many different random walks on the lattice. Just as Brownian motion can be viewed as 
a purely random continuous curve, the Brownian map seems to be the right model for a purely random surface. 
Note that the Brownian map $ ({\bf m}_\infty, D)$ is almost surely homeomorphic to the sphere $\mathbb{S}^2$ \cite{LGP} (see also \cite{Mi-sphere}),
even though its Hausdorff dimension is equal to $4$. However, 
there is no canonical homeomorphism and we cannot a priori use the Brownian map to obtain a ``canonical'' random metric on the sphere as suggested
by preceding remarks. 

The main goal of the present work is to present a general theorem of convergence toward the Brownian map
(Theorem \ref{mainresult} below), which was 
derived in the series of papers \cite{IM,AM,LGU}. We also give a detailed construction of the Brownian map, by showing
that this random compact metric space can be obtained by gluing certain pairs of points in another famous 
probabilistic model, the Brownian continuum random tree or CRT, which was introduced and studied by Aldous \cite{Al1,Al3}.
Note that the CRT itself is a universal scaling limit of discrete trees, in a sense analogous to (\ref{convtri}). The best way to
understand the construction of the Brownian map from the CRT is to start from the discrete bijections
that exist between various classes of planar maps and corresponding classes of  labeled trees. For this reason, we start
in Section \ref{secdiscrete} by explaining these bijections in the particular case of quadrangulations (there exist similar bijections
for triangulations, or more generally for $p$-angulations, but their description is more complicated).  In Section \ref{secBM}, after stating the main theorem of convergence
that extends (\ref{convtri}), we give a precise definition of the CRT as the tree coded by a normalized Brownian excursion, and
we construct the Brownian map via the above-mentioned gluing procedure. Note that, although the discrete bijections between planar
maps and trees are very far from being sufficient to get results such as (\ref{convtri}), they provide useful insight into the construction
of the Brownian map. In Section \ref{secgeodesic}, we describe the detailed results that are known about geodesics in the Brownian map.
Section \ref{secBP} is devoted to the Brownian plane, which can be viewed as an infinite-volume version of the Brownian map. The Brownian plane shares many properties of the Brownian map, but it enjoys an additional scale invariance property, which makes it
more suitable for explicit calculations (we briefly present some recent results from \cite{CLG2}). One can also view the Brownian plane as a continous analog of the infinite random 
lattices known as the uniform infinite planar triangulation and quadrangulation, which have been studied extensively in the
recent years. Finally, Section \ref{secopen} discusses what is probably the most important open problem
in the area, namely finding a canonical construction of the Brownian map in terms of a random metric on the sphere.
A related discussion can be found in Benjamini \cite{Be}.

To complete this introduction, let us mention that, although our main motivation for the following results came from
probability theory, there are important connections with several other areas of mathematics and physics. In particular, from the point of view
of combinatorics, one can derive information about properties of ``typical'' large planar maps from the known
results about the Brownian map. For instance, the fact that the Brownian map is homeomorphic to the sphere implies that
for a large planar map chosen uniformly at random in a given class (e.g. a large triangulation) there cannot exist a cycle whose
length is small in comparison with the diameter of the graph, such that both regions separated by this cycle have a 
macroscopic size (see \cite{LGP}).  Similarly, the statements about geodesics in the Brownian map show that, in a typical large planar map,
the geodesic between two vertices chosen at random is ``macroscopically unique'', meaning that the distance between two geodesics 
will be small in comparison with the diameter of the graph. Another strong motivation for this work came
from physics, and particularly from the use of large random planar maps as models of random geometry in two-dimensional
quantum gravity (see the book \cite{ADJ}). In this connection, we mention the important work of Bouttier and Guitter (see 
in particular \cite{BG1,BG2,BG3}) which motivated part of our results. It is worth mentioning that a different mathematical approach to
two-dimensional quantum gravity relying on the Gaussian free field has been given by Duplantier and Sheffield \cite{DS}. 
The construction of \cite{DS} seems quite different from our perspective, but 
the paper \cite{She} gives a number of conjectures that would relate large random planar maps and the Brownian map 
to the Gaussian free field approach. The very recent work of Miller and Sheffield \cite{MS} seems  promising in this
respect.  

\smallskip

\noindent{\it Acknowledgements.} I thank Nicolas Curien and Igor Kortchemski for their help with figures. 

\section{Planar maps and bijections with trees}
\label{secdiscrete}

In this section, we recall the basic definitions concerning planar maps and we explain how planar maps can be coded by trees, in the particular
case of quadrangulations. This coding is an important ingredient
of the proofs, and it also helps to understand the definition of the Brownian map that will
be given below.

\subsection{Planar maps}

Let us start with the definition of a planar map.

\begin{definition}
\label{def-planarmap}
A planar map is a proper (without edge-crossing) embedding of a finite connected graph in the 
two-dimensional sphere $\mathbb{S}^2$. A rooted planar map is a planar map given with a distinguished oriented
edge, which is called the root edge and whose tail vertex is called the root vertex. Two rooted planar maps
are identified if they correspond via an orientation-preserving homeomorphism of $\mathbb{S}^2$.
\end{definition}

We in fact allow loops and multiple edges in our graphs. Following the classical terminology found in combinatorics, the
word ``graph'' should be replaced by ``multigraph'' in the preceding definition. The faces of  a planar map
are the connected components of the complement of the union of edges, and the degree of a face counts
the number of edge sides that are incident to this face (in particular it may happen that both sides of an
edge are incident to the same face, and then this edge is counted twice in the degree of the face).

\medskip
\begin{figure}
\begin{center}
\ifx\JPicScale\undefined\def\JPicScale{.6}\fi
\unitlength \JPicScale mm
\begin{picture}(111,80)(0,0)
\linethickness{0.6mm}
\put(16,25){\line(0,1){30}}
\linethickness{0.6mm}
\put(16,55){\line(1,0){30}}
\linethickness{0.6mm}
\put(16,25){\line(1,0){30}}
\linethickness{0.6mm}
\put(46,25){\line(0,1){30}}
\linethickness{0.6mm}
\qbezier(16,55)(21.78,62.5)(37.88,66.25)
\qbezier(37.88,66.25)(53.97,70)(56,70)
\linethickness{0.6mm}
\qbezier(56,70)(58.66,57.34)(52.88,41.88)
\qbezier(52.88,41.88)(47.09,26.41)(46,25)
\linethickness{0.6mm}
\qbezier(56,70)(63.19,63.28)(64.75,40.62)
\qbezier(64.75,40.62)(66.31,17.97)(66,15)
\linethickness{0.6mm}
\qbezier(16,25)(38.81,11.25)(52.25,12.5)
\qbezier(52.25,12.5)(65.69,13.75)(66,15)
\linethickness{0.6mm}
\qbezier(56,70)(84.44,58.59)(88.5,51.88)
\qbezier(88.5,51.88)(92.56,45.16)(91,45)
\linethickness{0.6mm}
\put(91,25){\line(0,1){20}}
\linethickness{0.6mm}
\qbezier(66,15)(74.12,15.62)(82.25,20)
\qbezier(82.25,20)(90.38,24.38)(91,25)
\linethickness{0.6mm}
\qbezier(56,70)(63.97,77.03)(72.88,75.62)
\qbezier(72.88,75.62)(81.78,74.22)(86,70)
\qbezier(86,70)(89.59,65.16)(90.38,55.62)
\qbezier(90.38,55.62)(91.16,46.09)(91,45)
\linethickness{0.6mm}
\qbezier(56,70)(64.12,72.97)(72.25,69.38)
\qbezier(72.25,69.38)(80.38,65.78)(81,65)
\linethickness{0.6mm}
\qbezier(91,50)(107.09,34.53)(109.12,21.88)
\qbezier(109.12,21.88)(111.16,9.22)(86,5)
\qbezier(86,5)(59.75,1.25)(44.75,1.25)
\qbezier(44.75,1.25)(29.75,1.25)(26,5)
\qbezier(26,5)(22.09,9.53)(19.12,16.88)
\qbezier(19.12,16.88)(16.16,24.22)(16,25)
\linethickness{0.6mm}
\put(16.5,24.5){\circle*{5}}

\linethickness{0.6mm}
\put(16,54.5){\circle*{5}}

\linethickness{0.6mm}
\put(45.5,54.5){\circle*{5}}

\linethickness{0.6mm}
\put(46,25){\circle*{5}}

\linethickness{0.6mm}
\put(81,65){\circle*{5}}

\linethickness{0.6mm}
\put(56.5,70){\circle*{5}}

\linethickness{0.6mm}
\put(91,50){\circle*{5}}

\linethickness{0.6mm}
\put(91,25){\circle*{5}}

\linethickness{0.6mm}
\put(66,15){\circle*{5}}

\linethickness{0.7mm}
\multiput(11,32)(0.12,-0.14){42}{\line(0,-1){0.14}}
\linethickness{0.7mm}
\multiput(16,26)(0.12,0.14){42}{\line(0,1){0.14}}
\put(12,64){\makebox(0,0)[cc]{{\small root}}}

\put(10,60){\makebox(0,0)[cc]{{\small vertex}}}

\put(23,42){\makebox(0,0)[cc]{{\small root}}}

\put(23,38){\makebox(0,0)[cc]{{\small edge}}}

\end{picture}\end{center}

\caption{A rooted quadrangulation with $7$ faces}
\end{figure}
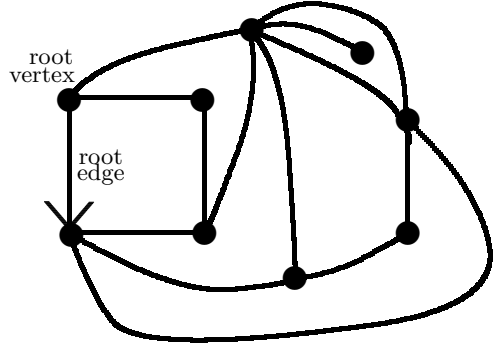

For any integer $p\geq 3$, a planar map is called a $p$-angulation 
(a triangulation if $p=3$, a quadrangulation if $p=4$)
if all its faces have degree $p$. Fig.~1 shows an example
of a quadrangulation with $7$ faces. For every integer $n\geq 1$, we write $\mathcal{M}^p_n$
for the set of all rooted $p$-angulations with $n$ faces. When $p$ is odd, the set $\mathcal{M}^p_n$ is empty if (and only if) $n$ is odd. So when we deal
with odd values of $p$ (in particular with triangulations), we will always assume that $n$ is even.

Thanks to the identification in Definition \ref{def-planarmap}, the sets $\mathcal{M}^p_n$ are finite. 
Enumeration results for these sets were 
obtained by Tutte (see in particular \cite{Tu}) in a series of important papers motivated by the four-color theorem. 

We will use the notation $d^M_{\rm gr}$ for the graph distance on the vertex set $V(M)$ of a planar map $M$.

\smallskip
\subsection{Bijections with trees}

In this section we explain the Cori-Vauquelin-Schaeffer bijection between
rooted quadrangulations and well-labeled trees \cite{CS,CV}.
We restrict ourselves to the special case of quadrangulations for the sake of simplicity, but
similar bijections exist for $p$-angulations for any $p$ (see in particular \cite{BDG}). 

First recall that a plane tree $\tau$ is a (finite) rooted ordered tree: A way of specifying a plane tree 
is to represent each of its vertices by a finite word made of positive integers, in such a way that 
the empty word $\varnothing$ corresponds to
the root or ancestor of the tree, and that, for instance, the word $13$ corresponds to the third child 
of the first child of the root (see the left side of Fig.~2). 

Then a well-labeled tree is a plane tree $\tau$, with vertex set $V(\tau)$, whose vertices $v$ are assigned 
labels $(\ell_v)_{v\in V(\tau)}$, in such a way that the following properties hold. First, the label of the root is $1$, then
the label $\ell_v$ of any vertex $v$ is a positive integer, and finally $|\ell_v -\ell_{v'}|\leq 1$ if 
$v$ and $v'$ are adjacent vertices.

With any well-labeled tree $(\tau,(\ell_v)_{v\in V(\tau)})$ with $n$ edges, we can associate a 
rooted quadrangulation $Q$ with $n$ faces via the following construction. 
Suppose that the tree is embedded in the plane in the (obvious) way as suggested in the 
left part of Fig.~2 (in particular the successive children of a vertex appear from left to right). 
Then the vertex set of $Q$ will be the union of the vertex set 
of $\tau$ and of an extra vertex, and we now explain how to construct the edges of $Q$.
First recall that 
a corner of the tree $\tau$ is an angular sector between two successive edges of $\tau$ around a given vertex.
The set of all corners of $\tau$ is given a cyclic ordering by moving clockwise around the tree.
To construct the edges of $Q$, we
first add an extra vertex $\partial$ outside the tree, and we connect each corner of the 
tree $\tau$ with label $1$ to the vertex $\partial$ by an edge starting from this corner. 
Then every corner of $\tau$ with label $k\geq2$ is connected by an edge to the next corner (in cyclic ordering)
with label $k-1$. The construction can be made in a unique way so
that edges do not cross and do not cross the edges of the tree. The resulting collection of edges  forms a quadrangulation $Q$
whose vertex set is $V(Q)=V(\tau)\cup\{\partial\}$. This quadrangulation is rooted at the edge
connecting the first corner of the root of $\tau$ to $\partial$, which is oriented so that $\partial$
is the root vertex. See Fig.~2 for an example.

\begin{figure}

\ifx\JPicScale\undefined\def\JPicScale{.9}\fi
\unitlength \JPicScale mm
\begin{picture}(120,80)(28,0)

\newsavebox{\CoriS}
\savebox{\CoriS}(100,70)[bl]{

\linethickness{0.3mm}
\put(80,18){\circle*{2}}

\linethickness{0.3mm}
\put(80,29){\circle*{2}}

\linethickness{0.3mm}
\put(89,36){\circle*{2}}

\linethickness{0.3mm}
\put(97,45){\circle*{2}}

\linethickness{0.3mm}
\put(71,36){\circle*{2}}

\linethickness{0.3mm}
\put(78,45){\circle*{2}}

\linethickness{0.3mm}
\put(63,45){\circle*{2}}

\linethickness{0.3mm}
\put(54,55){\circle*{2}}

\linethickness{0.2mm}
\put(80,18){\line(0,1){11}}
\linethickness{0.2mm}
\multiput(80,29)(0.16,0.12){58}{\line(1,0){0.16}}
\linethickness{0.2mm}
\multiput(89,36)(0.12,0.13){67}{\line(0,1){0.13}}
\linethickness{0.2mm}
\multiput(71,36)(0.16,-0.12){58}{\line(1,0){0.16}}
\linethickness{0.2mm}
\multiput(63,45)(0.12,-0.13){67}{\line(0,-1){0.13}}
\linethickness{0.2mm}
\multiput(54,55)(0.12,-0.13){75}{\line(0,-1){0.13}}
\linethickness{0.2mm}
\multiput(71,36)(0.12,0.16){58}{\line(0,1){0.16}}

\linethickness{0.3mm}

\put(76,16){\makebox(0,0)[cc]{\fbox{$\mathbf{1}$}}}
\put(82,16){\makebox(0,0)[cc]{$\scriptstyle{\varnothing}$}}

\put(76,27){\makebox(0,0)[cc]{\fbox{$\mathbf{2}$}}}
\put(82,28){\makebox(0,0)[cc]{$\scriptstyle{1}$}}

\put(67,35){\makebox(0,0)[cc]{\fbox{$\mathbf{3}$}}}
\put(74,37){\makebox(0,0)[cc]{$\scriptstyle{11}$}}

\put(78,49){\makebox(0,0)[cc]{\fbox{$\mathbf{3}$}}}
\put(81,43){\makebox(0,0)[cc]{$\scriptstyle{112}$}}

\put(59,43){\makebox(0,0)[cc]{\fbox{$\mathbf{2}$}}}
\put(67,46){\makebox(0,0)[cc]{$\scriptstyle{111}$}}

\put(50,54){\makebox(0,0)[cc]{\fbox{$\mathbf{1}$}}}
\put(58,56){\makebox(0,0)[cc]{$\scriptstyle{1111}$}}

\put(85,38){\makebox(0,0)[cc]{\fbox{$\mathbf{1}$}}}
\put(90,33){\makebox(0,0)[cc]{$\scriptstyle{12}$}}

\put(93,47){\makebox(0,0)[cc]{\fbox{$\mathbf{2}$}}}
\put(99,42){\makebox(0,0)[cc]{$\scriptstyle{121}$}}

}

\newsavebox{\Cori}
\savebox{\Cori}(100,70)[bl]{

\linethickness{0.3mm}
\put(80,18){\circle*{2}}

\linethickness{0.3mm}
\put(80,29){\circle*{2}}

\linethickness{0.3mm}
\put(89,36){\circle*{2}}

\linethickness{0.3mm}
\put(97,45){\circle*{2}}

\linethickness{0.3mm}
\put(71,36){\circle*{2}}

\linethickness{0.3mm}
\put(78,45){\circle*{2}}

\linethickness{0.3mm}
\put(63,45){\circle*{2}}

\linethickness{0.3mm}
\put(54,55){\circle*{2}}

\linethickness{0.2mm}
\put(80,18){\line(0,1){11}}
\linethickness{0.2mm}
\multiput(80,29)(0.16,0.12){58}{\line(1,0){0.16}}
\linethickness{0.2mm}
\multiput(89,36)(0.12,0.13){67}{\line(0,1){0.13}}
\linethickness{0.2mm}
\multiput(71,36)(0.16,-0.12){58}{\line(1,0){0.16}}
\linethickness{0.2mm}
\multiput(63,45)(0.12,-0.13){67}{\line(0,-1){0.13}}
\linethickness{0.3mm}
\put(113,35){\circle*{2}}

\linethickness{0.6mm}
\qbezier(89,36)(95.41,37.88)(96.38,41.25)
\qbezier(96.38,41.25)(97.34,44.62)(97,45)
\linethickness{0.6mm}
\qbezier(89,36)(103.62,74.66)(108.5,57.62)
\qbezier(108.5,57.62)(113.38,40.59)(113,35)
\linethickness{0.6mm}
\qbezier(63,45)(86.28,73.97)(88.38,57.38)
\qbezier(88.38,57.38)(90.47,40.78)(89,36)
\linethickness{0.2mm}
\multiput(71,36)(0.12,0.16){58}{\line(0,1){0.16}}
\linethickness{0.6mm}
\qbezier(71,36)(80.75,64.44)(80.75,49)
\qbezier(80.75,49)(80.75,33.56)(80,29)
\linethickness{0.2mm}
\multiput(54,55)(0.12,-0.13){75}{\line(0,-1){0.13}}
\linethickness{0.6mm}
\qbezier(54,55)(57.03,-7.28)(83.12,10.12)
\qbezier(83.12,10.12)(109.22,27.53)(113,35)
\linethickness{0.6mm}
\qbezier(54,55)(62.19,24.5)(70.75,25.5)
\qbezier(70.75,25.5)(79.31,26.5)(80,29)
\linethickness{0.6mm}
\qbezier(54,55)(57.84,47.81)(60.38,46.25)
\qbezier(60.38,46.25)(62.91,44.69)(63,45)
\linethickness{0.6mm}
\qbezier(71,36)(65.25,39.19)(64,42)
\qbezier(64,42)(62.75,44.81)(63,45)
\linethickness{0.6mm}
\qbezier(113,35)(99.78,26.06)(90.12,22)
\qbezier(90.12,22)(80.47,17.94)(80,18)
\linethickness{0.6mm}
\qbezier(89,36)(97.72,33.97)(105.12,34.38)
\qbezier(105.12,34.38)(112.53,34.78)(113,35)
\linethickness{0.6mm}
\qbezier(80,29)(83.28,23.72)(81.88,20.88)
\qbezier(81.88,20.88)(80.47,18.03)(80,18)
\linethickness{0.6mm}
\qbezier(71,36)(76.16,34.25)(78.12,31.75)
\qbezier(78.12,31.75)(80.09,29.25)(80,29)
\put(116,34){\makebox(0,0)[cc]{$\partial$}}

\put(78,16){\makebox(0,0)[cc]{$\mathbf{1}$}}

\put(77,29){\makebox(0,0)[cc]{$\mathbf{2}$}}

\put(69,34){\makebox(0,0)[cc]{$\mathbf{3}$}}

\linethickness{0.6mm}
\qbezier(78,45)(78.12,37.44)(79,33.25)
\qbezier(79,33.25)(79.88,29.06)(80,29)
\put(78,48){\makebox(0,0)[cc]{$\mathbf{3}$}}

\put(62,42){\makebox(0,0)[cc]{$\mathbf{2}$}}

\put(51,55){\makebox(0,0)[cc]{$\mathbf{1}$}}

\linethickness{0.6mm}
\qbezier(80,29)(83.84,34.69)(86.38,35.5)
\qbezier(86.38,35.5)(88.91,36.31)(89,36)
\put(87,38){\makebox(0,0)[cc]{$\mathbf{1}$}}

\put(96,48){\makebox(0,0)[cc]{$\mathbf{2}$}}
\linethickness{0.6mm}
\put(97,28){\line(0,-1){3}}
\put(97,25){\line(1,0){3}}

}

\put(-15,0){\usebox{\CoriS}}
\put(48,0){\usebox{\Cori}}
\end{picture}

\caption{The Cori-Vauquelin-Schaeffer bijection. On the left side, a well-labeled tree (the framed numbers are the 
labels assigned to the vertices).  On the right side, the edges of the associated quadrangulation $Q$ appear in thick curves. }
\end{figure}
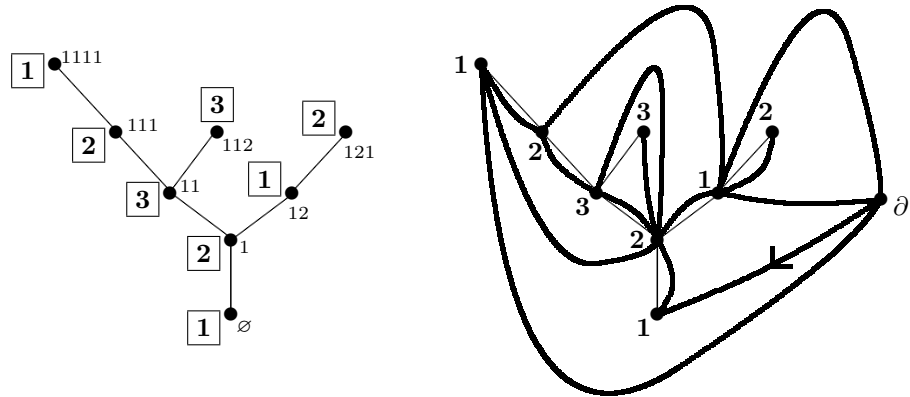

The previous construction yields a bijection from the set of all well-labeled trees with a fixed number $n$ of edges
onto the set $\mathcal{M}^4_n$ of all rooted quadrangulations with $n$ faces. This bijection is called the 
Cori-Vauquelin-Schaeffer bijection (the CVS bijection in short). Furthermore, the following
important additional property holds. If $(\tau,(\ell_v)_{v\in V(\tau)})$ is a well-labeled tree and $Q$ is the associated quadrangulation
defined as above, then, for every $v\in V(Q)\backslash\{\partial\}$,
\begin{equation}
\label{dist-root}
d^Q_{\rm gr}(\partial, v) = \ell_v.
\end{equation}
In other words, distances from the root vertex in the quadrangulation $Q$ are given by labels on the tree associated with $Q$ via the 
CVS bijection. There is no similar expression for $d^Q_{\rm gr}(v,v')$ when $v$ and $v'$ are two vertices other than $\partial$, but
the following upper bound turns out to be important for our purposes. For every $v,v'\in V(Q)\backslash\{\partial\}$,
\begin{equation}\label{bounddist}
d^Q_{\rm gr}(v,v')\leq \ell_v + \ell_{v'} - 2 \max\Big( \min_{w\in[v,v']}\ell_w , \min_{w\in[v',v]}\ell_w\Big) + 2,
\end{equation}
where $[v,v']$ stands for the set of all vertices visited when going from $v$ to $v'$ in clockwise order around the tree
(for instance, in the tree of Fig.~2, if $v= 111$ and $v'=12$, $[v,v']=\{111,11,112,1,12\}$
and $[v',v]=\{12,1,\varnothing,11,111\}$). The bound (\ref{bounddist})
is easily derived by the following argument. Consider any corner of $v$. Via the CVS bijection, this corner is connected
by an edge of $Q$
to a corner of another vertex $v_1$ with label $\ell_v-1$, and then this corner of $v_1$ is connected by an edge of $Q$
to a corner of a second vertex $v_2$ with label $\ell_v-2$. Recalling that labels correspond to distances 
from $\partial$, we get a geodesic path $\gamma=(v,v_1,v_2,\ldots)$ from $v$ to $\partial$. We may construct a similar geodesic path
$\gamma'$ from $v'$ to $\partial$ and observe that the geodesic paths $\gamma$ and $\gamma'$ will 
eventually merge. Considering the path
from $v$ to $v'$ obtained by concatenating the parts of $\gamma$ and $\gamma'$ before their merging point (and choosing the 
initial corners in an optimal way) easily leads
to the bound (\ref{bounddist}). 

The basic underlying idea of our construction of the Brownian map in the next section
is to use a continuous analog of the CVS bijection. In this analog, the role of
plane trees will be played by Aldous' continuum random tree, which is known
to be the (universal) scaling limit of many different classes of random discrete trees.
We refer to \cite{LGM} for a discussion of scaling limits of labeled plane trees.

\section{The Brownian map}
\label{secBM}

\subsection{The Gromov-Hausdorff distance}
\label{GHsubsec}

Let
$(E_1,d_1)$ and $(E_2,d_2)$ be two compact metric spaces. 
The Gromov-Hausdorff distance between $(E_1,d_1)$
and $(E_2,d_2)$ is
$$d_{GH}(E_1,E_2)=\inf\Big(d_{\rm Haus}(\varphi_1(E_1),\varphi_2(E_2))\Big),$$ 
where the infimum is over all isometric embeddings $\varphi_1:E_1\longrightarrow E$ and
$\varphi_2:E_2\longrightarrow E$ of $E_1$ and $E_2$ into the same metric
space $(E,d)$, and $d_{\rm Haus}$ stands for the usual Hausdorff distance 
between compact subsets of $E$. If $\mathbb{K}$ denotes the space of all
isometry classes of compact metric spaces, then $d_{GH}$ is a distance
on $\mathbb{K}$, and moreover the metric space $(\mathbb{K},d_{GH})$ is Polish, that is, separable
and complete. We refer to Chapter 7 of Burago, Burago and Ivanov \cite{BBI} for a thorough
discussion of the Gromov-Hausdorff distance.

\subsection{The main theorem}

Recall our notation $\mathcal{M}^p_n$ for the space of all rooted
$p$-angulations with $n$ faces. Note that when $p$ is odd we consider
only even values of $n$. With every $M\in \mathcal{M}^p_n$, we can associate  the
metric space $(V(M), d^M_{\rm gr})$.

\begin{theorem}[\cite{LGU}]
\label{mainresult}
Suppose that either $p=3$ or $p\geq 4$ is an even integer, and set
$$c_p:= \Big(\frac{9}{p(p-2)}\Big)^{1/4}$$
if $p$ is even, and 
$$c_3:= 6^{1/4}.$$
For every integer $n\geq1$ (for every even integer $n\geq 2$ if $p=3$), let $M_n$ be uniformly distributed over $\mathcal{M}^p_n$.
There exists a random compact metric space $(\mathbf{m}_\infty,D)$ called the Brownian map,
which does not depend on $p$, such that
$$(V(M_n), c_p\,n^{-1/4}d^{M_n}_{\rm gr}) \build{\longrightarrow}_{n\to\infty}^{\rm(d)} (\mathbf{m}_\infty, D)$$
where the convergence holds in distribution in the space $(\mathbb{K},d_{GH})$. 
\end{theorem}

The role of the scaling constants $c_p$ is only to ensure that the limit does not depend on $p$. 
The convergence in distribution of the theorem may be rephrased by saying that 
one can construct the full sequence $(M_n)$ in such a way that 
the associated metric spaces $(V(M_n), c_p\,n^{-1/4}d^{M_n}_{\rm gr})$ converge to $(\mathbf{m}_\infty, D)$
for the Gromov-Hausdorff distance, outside a set of zero probability. The name Brownian map is
due to Marckert and Mokkadem \cite{MaMo}, who proved a weak form of the theorem when $p=4$.
As was already pointed in the introduction,  the case $p=4$ of the theorem has been derived independently
by Miermont \cite{Mi-quad} using a different approach. 

\vspace{-1mm}

\begin{figure}[h]
\begin{center}
\includegraphics[width=12cm,height=8cm]{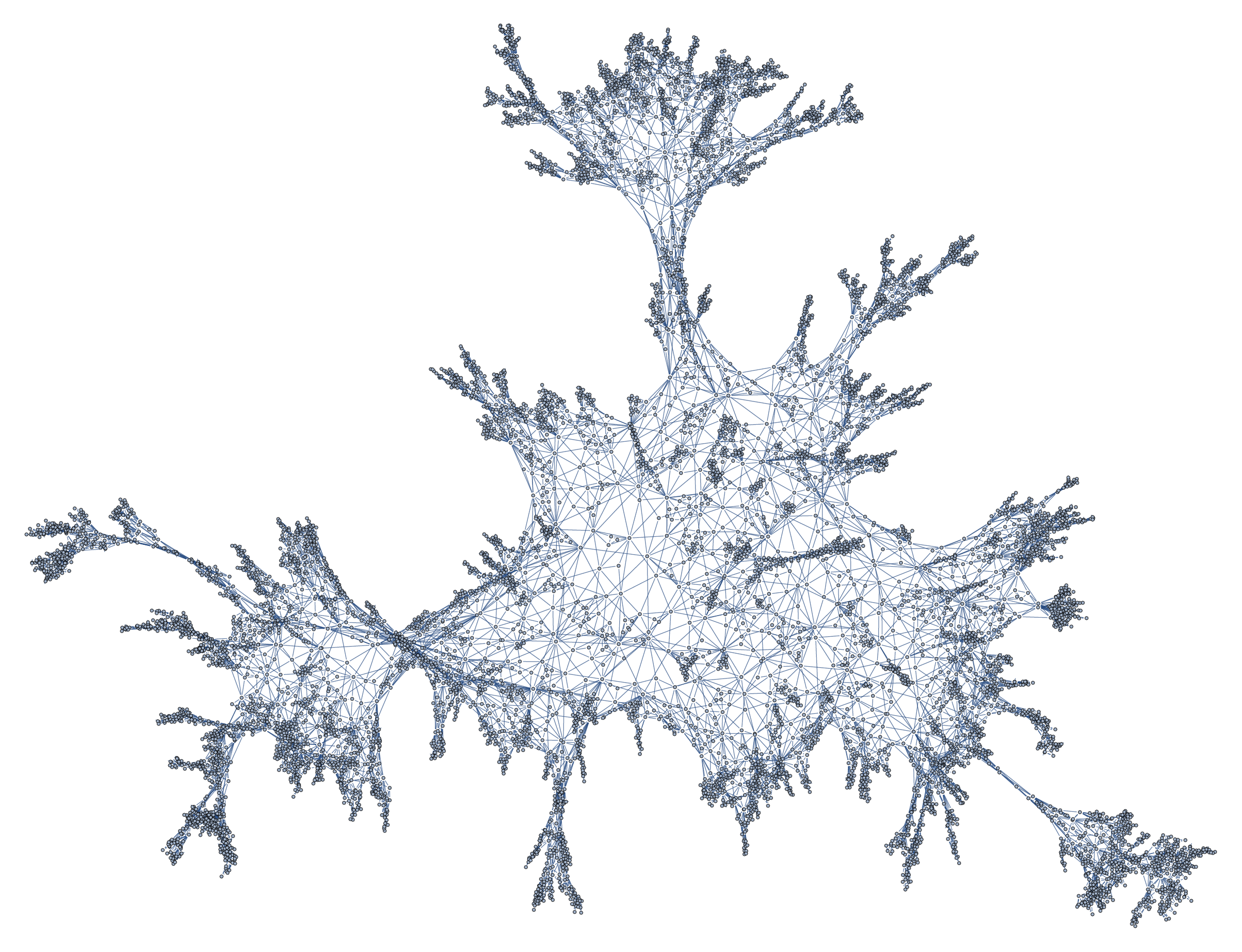}
\end{center}
\vspace{-2mm}
\caption{Simulation of a large triangulation of the sphere. Here only the graph structure is represented in three
dimensions.}
\end{figure}

\noindent{\bf Remark.} It is very likely that the convergence of the theorem also holds for odd values of $p\geq 5$, though additional
technical difficulties appear in that case.
Similarly, it is expected that
a version of the convergence holds for Boltzmann distributed 
random planar maps, such that the probability of a given planar map $M$ (with a fixed 
number $n$ of vertices) will be proportional to
$$\prod_{f\;{\rm face\;of\;}M} w_{{\rm degree}(f)}$$
where $(w_k)_{k\geq 1}$ is a suitable sequence of weights. In the bipartite case where
$w_k=0$ when $k$ is odd, a version of the convergence of the theorem
holds for such random planar maps \cite{LGU} under appropriate assumptions. Other recent extensions of the
theorem have been mentioned in the introduction above.

\smallskip

In the next two subsections, we will present a precise construction of the Brownian map as a quotient space of
another (well-known) random compact metric space, which is the Brownian continuum random tree
or CRT (see Aldous \cite{Al1,Al3}). 

\subsection{The Brownian continuum random tree}

We first recall the notion of an $\mathbb{R}$-tree.

\begin{definition}
\label{defRTree}
A metric space $(\mathcal{T},d)$ is an $\mathbb{R}$-tree if the following two
properties hold for every $a,b\in \mathcal{T}$.
\begin{description}
\item{\rm(a)} There is a unique
isometric map
$f_{a,b}$ from $[0,d(a,b)]$ into $\mathcal{T}$ such
that $f_{a,b}(0)=a$ and $f_{a,b}(
d(a,b))=b$.
\item{\rm(b)} If $q$ is a continuous injective map from $[0,1]$ into
$\mathcal{T}$, such that $q(0)=a$ and $q(1)=b$, we have
$$q([0,1])=f_{a,b}([0,d(a,b)]).$$
\end{description}
\noindent A rooted $\mathbb{R}$-tree is an $\mathbb{R}$-tree $(\mathcal{T},d)$
with a distinguished vertex $\rho$ called the root.
\end{definition}

We will be interested mainly  in compact $\mathbb{R}$-trees.
Informally, one should think of a compact $\mathbb{R}$-tree as a 
connected union of line segments
in the plane with  no loops, which is equipped with the appropriate
(intrinsic) metric. For any two
points $a$ and $b$ in the tree, there is a unique arc going from $a$ to $b$
in the tree, which is isometric to a line segment. 

Rooted $\mathbb{R}$-trees can be coded by contour functions, in a way very similar
to the 
well-known coding of plane trees by Dyck paths. Let 
$g:[0,1]\longrightarrow \mathbb{R}_+$ be a nonnegative continuous function such
that $g(0)=g(1)=0$ and, for every $s,t\in[0,1]$, set
$$m_g(s,t):=\inf_{r\in[s\wedge t,s\vee t]}g(r),$$
and
$$d_g(s,t):=g(s)+g(t)-2m_g(s,t).$$
It is easy to verify that $d_g$ is a pseudo-metric 
on $[0,1]$. As usual, we introduce the equivalence
relation
$s\sim_g t$ if and only if $d_g(s,t)=0$ (or equivalently if and only if $g(s)=g(t)=m_g(s,t)$).
The function $d_g$ induces a distance on the quotient space $\mathcal{T}_g
:=[0,1]\,/\!\sim_g$, and we keep the
notation $d_g$ for this distance. We also write $p_g$ for the
canonical projection from $[0,1]$ onto $\mathcal{T}_g$. Then it is not hard to verify that 
the metric space $(\mathcal{T}_g,d_g)$ is a compact $\mathbb{R}$-tree (see e.g.~\cite{DuL}), which by definition is rooted at
$\rho_g=p_g(0)=p_g(1)$. 
Furthermore the mapping $g\longrightarrow \mathcal{T}_g$
is continuous with respect to the Gromov-Hausdorff distance,
if the set of continuous functions $g$ is equipped with the supremum distance.

The preceding coding induces a cyclic ordering on the tree $\mathcal{T}_g$, and it will
be important for us to consider the corresponding intervals. By convention,
if $s,t\in[0,1]$ are such that $s>t$, we set $[s,t]:=[s,1]\cup[0,t]$. 
Then we note that,
for every $a,b\in \mathcal{T}_g$ with $a\not = b$, there exists a smallest interval $[s,t]$
such that $p_g(s)=a$ and $p_g(t)=b$, and we define $[a,b]:=p_g([s,t])$. Roughly speaking,
$[a,b]$ corresponds to the set of vertices that are visited when going
from $a$ to $b$ in ``clockwise order around the tree''. 

Let ${\bf e}=({\bf e}_t)_{0\leq t\leq 1}$ be a normalized Brownian excursion. Informally,
${\bf e}$ behaves like a linear Brownian started from $0$ at time $0$, which is conditioned
to stay positive over the time interval $(0,1)$ and to return to $0$ at time $1$ (of course these
conditionings require special care, see e.g. \cite[Chapter XII]{RY} for a rigorous definition and many
properties of the Brownian excursion).

\begin{definition}
\label{def-CRT}
The CRT is the random $\mathbb{R}$-tree $(\mathcal{T}_{\bf e},d_{\bf e})$ coded by the Brownian excursion ${\bf e}$.
\end{definition}

\begin{figure}[h]
\begin{center}
\includegraphics[width=12cm,height=8cm]{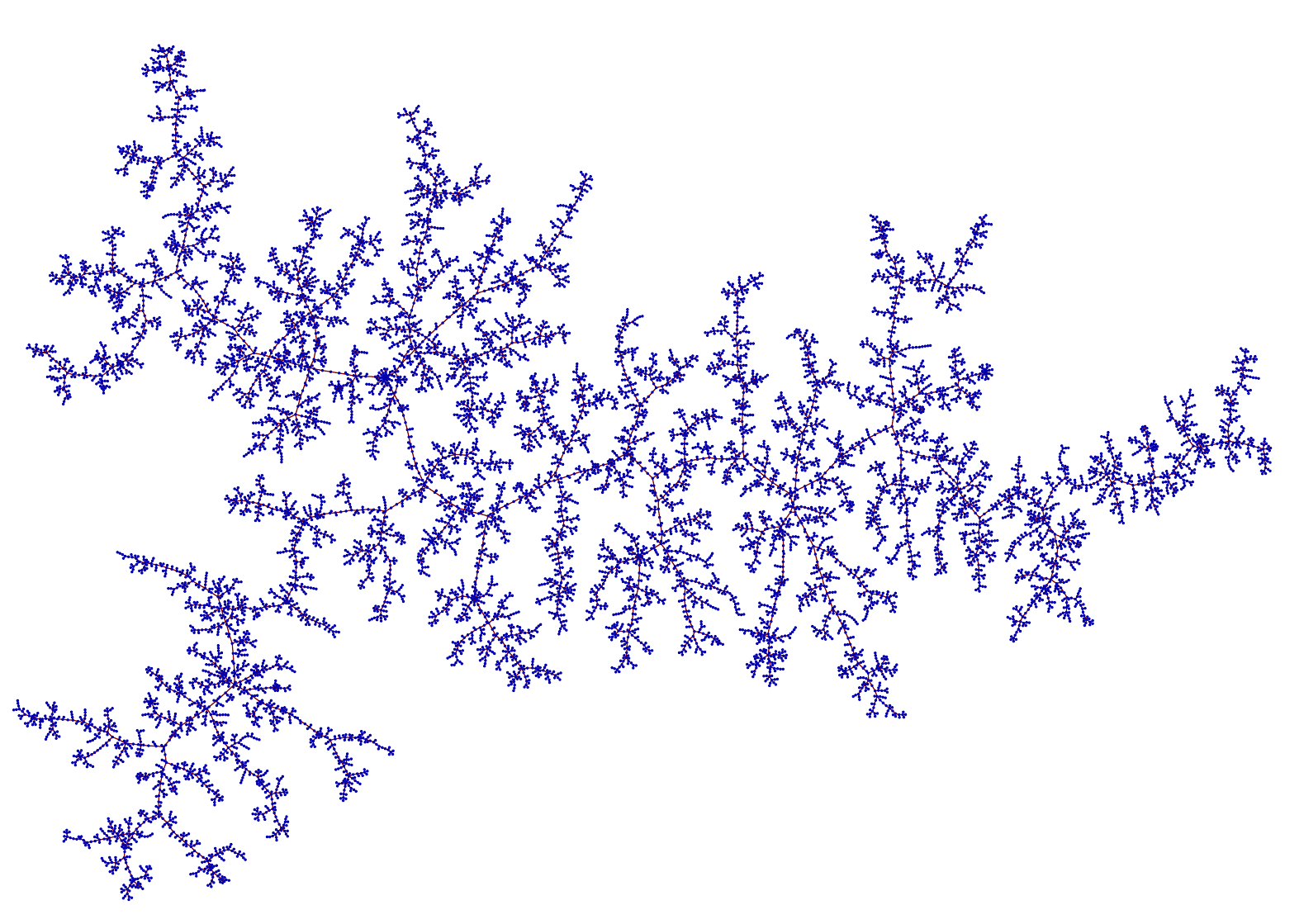}
\end{center}
\caption{A large discrete random tree, which is an approximation of the CRT $\mathcal{T}_{\mathbf{e}}$.}
\end{figure}

\subsection{Constructing the Brownian map} 
\label{consBM}
Analogously to the construction of
quadrangulations via the CVS bijection, we will need to introduce labels on $\mathbb{R}$-trees.
Consider first a deterministic $\mathbb{R}$-tree $(\mathcal{T},d)$, which is rooted at $\rho$. We define Brownian 
motion indexed by $\mathcal{T}$ as the centered Gaussian process
$(\mathcal{Z}_a)_{a\in\mathcal{T}}$ whose distribution is characterised by the properties $\mathcal{Z}_\rho=0$ and
$$E[(\mathcal{Z}_a-\mathcal{Z}_b)^2] = d(a,b)$$
for every $a,b\in \mathcal{T}$. For a general $\mathbb{R}$-tree, $(\mathcal{Z}_a)_{a\in\mathcal{T}}$ needs not have continuous paths.
However, if $\mathcal{T}=\mathcal{T}_g$ and the function $g$ is H\"older continuous, then it is easy to verify that
we can construct $(\mathcal{Z}_a)_{a\in\mathcal{T}_g}$ so that it has continuous sample paths. This applies 
in particular to (almost every realization of) the CRT $\mathcal{T}_{\bf e}$.

The building blocks of our construction are first the CRT $(\mathcal{T}_{\bf e},d_{\bf e})$
and then the process $(Z_a)_{a\in\mathcal{T}_{\bf e}}$ which, conditionally on the CRT, is 
distributed as Brownian motion indexed by $\mathcal{T}_{\bf e}$. Note that there seems to
be a technical difficulty here since we are considering a random process $Z$ indexed
by a random set $\mathcal{T}_{\bf e}$: A simple way out is first to construct a process $(Z_t)_{t\in[0,1]}$ as
the ``tip''  of the Brownian snake driven by ${\bf e}$ (see \cite{LZurich}), and then to 
observe that $Z_s=Z_t$ if $s\sim_{\bf e} t$, which allows one to view $Z$ as indexed by
$\mathcal{T}_{\bf e}$. 

We set, for every $a,b\in\mathcal{T}_{\bf e}$,
\begin{equation}
\label{Dzero}
D^\circ(a,b):=Z_a+Z_b-2\max\Big(\min_{c\in[a,b]}Z_c,\min_{c\in[b,a]}Z_c\Big).
\end{equation}
Note that we are using the notion of a tree interval which was introduced in the previous
subsection, and that $D^\circ(a,b)$ is a continuous analog of the right side of (\ref{bounddist}). 

The function $D^\circ$ is symmetric but does not satisfy the triangle inequality. We then consider
the largest symmetric function that is bounded above by $D^\circ$ and satisfies the
triangle inequality:
for every $a,b\in \mathcal{T}_\mathbf{e}$,
\begin{equation}
\label{formulaD}
D(a,b):= \inf\Big\{ \sum_{i=1}^k D^\circ(a_{i-1},a_i)\Big\},
\end{equation}
where the infimum is over all choices of the integer $k\geq 1$ and of the
elements $a=a_0,a_1,\ldots,a_{k-1},a_k=b$ of $\mathcal{T}_\mathbf{e}$.

Then $D$ is a pseudo-metric on $\mathcal{T}_{\bf e}$, and we set
$$a\approx b\quad\hbox{if and only if}\quad D(a,b)=0.$$

\begin{definition}
\label{def-Brownianmap}
The Brownian map is the quotient space ${\bf m}_\infty: = \mathcal{T}_{\bf e}\,/\!\approx$ equipped with the
distance induced by $D$.
\end{definition}

We will write $\Pi$ for the canonical projection from $\mathcal{T}_{\bf e}$ onto ${\bf m}_\infty$,
and we keep the notation $D$ for the induced distance on ${\bf m}_\infty$,
so that $D(\Pi(a),\Pi(b))=D(a,b)$ for every $a,b\in\mathcal{T}_{\bf e}$.

One can prove \cite{IM} that the property $D(a,b)=0$ holds if and only if $D^\circ(a,b)=0$, or
equivalently
\begin{equation}
\label{ident-vertex}
Z_a=Z_b=\max\Big(\min_{c\in[a,b]}Z_c,\min_{c\in[b,a]]}Z_c\Big).
\end{equation}
In other words, two vertices of the CRT are identified if they have the
same label $Z_a=Z_b$ and if one can go from $a$ to $b$ ``around the tree''
(in clockwise or couterclockwise order) encountering only vertices with label 
greater than or equal to $Z_a$. In a sense, not many pairs of vertices are
identified. Only leaves of the CRT may be identified (a leaf is a vertex $a$ of $\mathcal{T}_{\bf e}$
such that $\mathcal{T}_{\bf e}\backslash\{a\}$ remains connected), and the set of all vertices $a$ that are identified to another one has 
Hausdorff dimension $1$ (whereas the CRT had dimension $2$). Furthermore, there
is only a countable collection
of equivalence classes of $\approx$ that contain $3$ points, and
no equivalence class contains more than $3$ points. Still these identifications
drastically change the topology. 

\medskip
\noindent{\bf Remark.} As already mentioned above, the preceding construction of the Brownian map is 
analogous to the construction of quadrangulations via well-labeled trees
in the CVS bijection (and to the construction of more general planar maps
using other similar bijections). The pair $(\mathcal{T}_{\bf e},(Z_a)_{a\in\mathcal{T}_{{\bf e}}})$
is a kind of continuous analog of a well-labeled tree $(\tau,(\ell_v)_{v\in V(\tau)})$ 
(there is a minor difference because we do not impose a positivity condition
on the labels in the continuous setting, but we could have re-rooted $\mathcal{T}_{\bf e}$
at the vertex with minimal label and shifted all labels $Z_a$ so that they
become nonnegative, which would not have affected the preceding
construction). Still one may notice that no identification of vertices is
needed in the discrete setting. The reason why such identifications become
necessary in the continuous setting can be explained intuitively as follows.
In a large well-labeled tree, there will exist vertices $u$ and $v$, which are
at a ``macroscopic distance'' in the tree, and such that
$\ell_v=\ell_{u}-1$ and $\ell_w\geq \ell_u$ for every vertex $w\in[u,v]$, where
the interval $[u,v]$ is as in (\ref{bounddist}) (note that (\ref{ident-vertex}) is just a continuous version of
these properties): According to the rules of the CVS bijection, two such 
vertices are linked by an edge of the quadrangulation, and asymptotically,
recalling that the graph distance is rescaled by a factor tending to $0$, this
leads to an identification of two vertices in the tree. 

\subsection{Properties of the Brownian map}
In this subsection, we describe a few properties of the Brownian map, which show that this
random metric space behaves very differently from a smooth surface.

\begin{proposition}
\label{proper}
{\rm (i)} The Brownian map $({\bf m}_\infty, D)$ is a.s. homeomorphic to the sphere $\mathbb{S}^2$.
\par\noindent{\rm (ii)} The Hausdorff dimension of $({\bf m}_\infty, D)$ is a.s. equal to $4$.
\end{proposition}

Part (i) is proved in \cite{LGP}, and part (ii) was derived in \cite{IM}.
In view of Proposition \ref{proper}, one may look at the Brownian map as a sphere with a very
singular metric. Property (ii) is of course reminiscent of the fact that the Hausdorff dimension
of  a Brownian path in $\mathbb{R}^d$, $d\geq2$ is twice the dimension of a smooth curve.

We next turn to uniform estimates on the volume of balls. We first need to introduce 
the volume measure on ${\bf m}_\infty$, which in a sense is uniformly distributed over this
random space. We set ${\bf p}:=\Pi\circ p_{\bf e}$, which maps $[0,1]$ onto ${\bf m}_\infty$. The volume measure 
$\Lambda$ on
the Brownian map is the image of Lebesgue measure on $[0,1]$ under ${\bf p}$. We could also have 
introduced this volume measure via the convergence in Theorem \ref{mainresult}, 
as the limit of the normalized counting measure on the vertex set of $M_n$
(to make this rigorous, one should state the latter convergence in the sense of the 
Gromov-Hausdorff topology for measured metric spaces). 

If $x\in{\bf m}_\infty$ and $r>0$, we write $B_D(x,r)$ for the closed ball of radius $r$
centered at $x$ in $({\bf m}_\infty,D)$.

\begin{proposition}
\label{estimate-volume}
Let $\delta>0$. There exist two (random) positive constants $c_\delta$ and $C_\delta$
such that, for every $r\in(0,1]$ and every $x\in {\bf m}_\infty$,
$$c_\delta\,r^{4+\delta} \leq \Lambda(B_D(x,r))\leq C_\delta\,r^{4-\delta}.$$
\end{proposition}

This shows that in a way the Brownian map is ``very regular in its irregularity'': The volume
of any small ball of radius $r$ is approximately $r^4$. The upper bound in Proposition
\ref{estimate-volume} is proved in \cite{AM}, and the lower bound is easy from the
bound $D\leq D^\circ$ and the H\"older continuity properties of the process $Z$. 

Another way of quantifying the irregularity of the Brownian map is to look at
connected components of the complement of a ball. If this complement if not empty, it
will have typically infinitely many connected components. The following proposition from
\cite{LGcactus} gives a more precise estimate.

\begin{proposition}
\label{cactus}
Suppose that $U$ is a point of ${\bf m}_\infty$
chosen at random according to the volume measure $\Lambda$. Let $r>0$ and, 
for every $\varepsilon>0$, let $N^r_\varepsilon$ be the number of connected components
of $B_D(U,r)^c$ that intersect $B_D(U,r+\varepsilon)^c$. Then,
$$\varepsilon^3\, N^r_\varepsilon \build\longrightarrow_{\varepsilon \to0}^{} c_0\,\lambda_r,$$
where the convergence holds in probability,  $c_0>0$ is an explicit constant and $\lambda_r$ is the value at $r$ of the continuous density of the profile 
of distances from $U$, which is the random measure $\mathcal{I}$ on $\mathbb{R}_+$ defined by
$$\mathcal{I}(A)=\int \Lambda(dx)\, \mathbf{1}_A(D(U,x)),$$
for any Borel subset $A$ of $\mathbb{R}_+$. 
\end{proposition}

\section{Geodesics in the Brownian map}
\label{secgeodesic}

Although the Brownian map remains a mysterious object in many respects, there is detailed 
information about the structure of geodesics toward a typical point. In this section we describe
these results following \cite{AM}. We rely on the construction given in subsection \ref{consBM}. 

One can prove (see in particular \cite{LGW}) that there exists a unique vertex $\rho_*\in\mathcal{T}_{{\bf e}}$
such that
$$Z_{\rho_*}= \min_{c\in\mathcal{T}_{\bf e}} Z_c.$$
We set $x_*:=\Pi(\rho_*)$. In what follows we discuss geodesics toward $x_*$
in the Brownian map. It should be emphasized that $x_*$
is not a special point of the Brownian map: If $U$ is a point of the Brownian 
map chosen at random according to the volume measure, one can check \cite{AM} that
the random pointed compact metric spaces $({\bf m}_\infty, D, U)$ and $({\bf m}_\infty, D, x_*)$
have the same distribution, so that our results also give information about geodesics
toward the ``typical'' point $U$.

To simplify notation, we set
$$Z_*:=\min_{c\in\mathcal{T}_{\bf e}} Z_c.$$
We first observe that, for every $a\in\mathcal{T}_{{\bf e}}$,
\begin{equation}
\label{dist-root-BM}
D(x_*,\Pi(a)) = Z_a - Z_*.
\end{equation}
Indeed, it immediately follows from (\ref{Dzero}) that $D(x_*,\Pi(a))\leq D^\circ(\rho_*,a)= Z_a-Z_*$. On the other hand, 
by using
the easy bound
\begin{equation}
\label{easy-bd}
D^\circ(a,b)\geq |Z_a-Z_b|,
\end{equation}
which also follows from (\ref{Dzero}), we immediately get from (\ref{formulaD}) that
the reverse inequality $D(x_*,\Pi(a))=D(\rho_*,a)\geq Z_a-Z_*$ holds.

Let $x=\Pi(a)$, $a\in\mathcal{T}_{\bf e}$ be a point in the Brownian map. We can construct certain geodesics
from $x$ to $x_*$ in the following manner. Choose $t\in[0,1)$ such that $p_{\bf e}(t)= a$, 
and for every $r\in[0,Z_a-Z_*]$, set
$$\gamma_t(r):=\left\{
\begin{array}{ll}
\!\min\{s\in[t,1]: Z_{p_{\bf e}(s)}=Z_a-r\}&\hbox{if }\{s\in[t,1]: Z_{p_{\bf e}(s)}=Z_a-r\}\not =\varnothing,\\
\!\min\{s\in[0,t]: Z_{p_{\bf e}(s)}=Z_a-r\}&\hbox{otherwise}.
\end{array}
\right.
$$
Informally, $p_{\bf e}(\gamma_t(r))$ is the first vertex with label $Z_a -r$ that one encounters 
when exploring the tree clockwise  starting from $a$. Note that, when $a$
is not a leaf, there are more than one way of starting from $a$: This corresponds to
the different possible choices of $t$. Since $\rho_*$ is the unique
vertex with minimal label, it is clear that $p_{\bf e}(\gamma_t(Z_a-Z_*))=\rho_*$. 

\begin{lemma}
\label{simple-geodesic}
Let $a\in \mathcal{T}_{{\bf e}}$ and $t\in[0,1)$ such that $p_{\bf e}(t)=a$. Set $\Gamma_t(r)={\bf p}(\gamma_t(r))$
for every $r\in [0,Z_a-Z_*]$. Then $(\Gamma_t(r))_{0\leq r\leq Z_a-Z_*}$ is a geodesic
from $\Pi(a)$ to $x_*$ in $({\bf m}_\infty,D)$. Such geodesics are called simple geodesics.
\end{lemma}

The proof is easy. If $0\leq r\leq r'\leq Z_a-Z_*$, the bound $D\leq D^\circ$ and the definition (\ref{Dzero}) immediately shows that
$D(\Gamma_t(r),\Gamma_t(r'))\leq D^\circ(p_{\bf e}(\gamma_t(r)),p_{\bf e}(\gamma_t(r')))= r'-r$. On the other
hand, by (\ref{dist-root-BM}), $D(\Gamma_t(0),\Gamma_t(Z_a-Z_*))= D(\Pi(a),x_*)=Z_a-Z_*$, and the 
triangle inequality now gives $D(\Gamma_t(r),\Gamma_t(r'))=r'-r$ for every $0\leq r\leq r'\leq Z_a-Z_*$.

Simple geodesics are indeed analogs in our continuous setting of the discrete geodesics for
quadrangulations that were briefly discussed at the end of Section \ref{secdiscrete}. 
The following proposition \cite{AM} is the key result in our study of geodesics.

\begin{proposition}
\label{simple-geo}
All geodesics to $x_*$ in $({\bf m}_\infty,D)$ are simple geodesics.
\end{proposition}

This proposition makes it possible to classify all the geodesics to $x_*$. Indeed, it is 
easy to count simple geodesics. If we fix a point $x\in{\bf m}_\infty$, a simple geodesic $\gamma_t$
from $x$ to $x_*$ is obtained by choosing first  $a\in\mathcal{T}_{\bf e}$ such that $\Pi(a)=x$, and
then $t\in[0,1)$ such that $p_{\bf e}(t)=a$. The choice of $a$ is in fact irrelevant: If $\Pi^{-1}(x)$ is not 
a singleton, then the vertices $a$ in $\Pi^{-1}(x)$ must be leaves, then, for each such $a$, there is
a unique $t\in[0,1)$ with $p_{\bf e}(t)=a$, and one immediately verifies that the associated
simple geodesics coincide. On the other hand, if $x=\Pi(a)$ and $a$ is not a leaf of $\mathcal{T}_{\bf e}$, then
there are ($2$ or $3$) values of $t$ such that $p_{\bf e}(t)=a$, and the corresponding simple
geodesics are distinct. The preceding discussion leads to the following theorem.

\begin{theorem}
\label{structure-geodesic}
Let ${\rm Sk}(\mathcal{T}_{\bf e})$ stand for the set of all vertices of the CRT that are not leaves,
and ${\rm Skel}:=\Pi({\rm Sk}(\mathcal{T}_{\bf e}))$. Then the restriction of $\Pi$ to ${\rm Sk}(\mathcal{T}_{\bf e})$
is a homeomorphism, and the Hausdorff dimension of ${\rm Skel}$
is equal to $2$. Furthermore, a.s. for every $x\in {\bf m}_\infty$,
\begin{enumerate}
\item[$\bullet$] if $x\in {\bf m}_\infty \backslash {\rm Skel}$, there is a unique geodesic  from $x$ to $x_*$;
\item[$\bullet$] if $x\in {\rm Skel}$, the number of distinct geodesics from $x$ to $x_*$ is the
multiplicity of $x$ in ${\rm Skel}$, that is, the number of connected components of ${\rm Skel}\backslash\{x\}$.
This multiplicity is either $2$ or $3$. 
\end{enumerate}
\end{theorem}

The set ${\rm Skel}$, which is a dense subset of the Brownian map homeomorphic to
a non-compact real tree, thus appears as the cut-locus of the Brownian map with respect 
to the point $x_*$: A point $x$ belongs to ${\rm Skel}$ if and only if there are at least two distinct
geodesics from $x$ to $x_*$. Perhaps suprisingly, there is a strong analogy with classical
results of differential geometry that go back to Poincar\'e. For a smooth surface homeomorphic 
to the sphere, the cut-locus is also a tree, and the number of distinct geodesics from a 
point of the cut-locus is equal to its multiplicity. 

It is easy to verify that the set ${\rm Skel}$ has zero volume measure (in terms of Hausdorff dimension, 
it is already clear that ${\rm Skel}$ is a ``small'' subset of ${\bf m}_\infty$). A consequence of Theorem \ref{structure-geodesic}
is thus the fact that, if one picks independently two points $x$ and $x'$ according to the volume measure
on ${\bf m}_\infty$, then a.s. there is a unique geodesic between $x$ and $x'$
(see also Miermont \cite{Mi2}  for a related result with a different approach). From this and Theorem \ref{mainresult},
one can deduce a property of ``macroscopic uniqueness'' of discrete geodesics in large planar maps, which was
already mentioned in the introduction. See \cite{AM} for more details. 

Theorem \ref{structure-geodesic} shows that our construction of the Brownian map 
as a quotient space of the CRT has a strong geometric meaning 
(although it is certainly not the only possible construction). Indeed the set ${\rm Skel}$,
which is a homeomorphic image of the skeleton of the CRT in our construction, is
a  geometric object defined intrinsically in terms of the Brownian map, and thus does not
depend on the particular construction we have developed. 

\smallskip
\noindent{\bf Remark}. The version of Theorem \ref{structure-geodesic} given in \cite{AM} applies to any random 
compact metric space which appears as a Gromov-Hausdorff limit in distribution of rescaled
$2p$-angulations (the uniqueness of such a limit was not yet known). This description of
geodesics was then a key ingredient of the proof of Theorem \ref{mainresult} in \cite{LGU}

\smallskip

The following corollary gives a confluence property of geodesics, which easily follows 
from the fact that two simple geodesics will always merge before their endpoint.

\begin{corollary}
\label{confluence}
Let $\delta >0$. Then a.s. there exists $\varepsilon>0$ such that, whenever $x$ and $x'$
are two points of ${\bf m}_\infty$ with $D(x_*,x)\geq \delta$ and $D(x_*,x')\geq \delta$, then,
if $f$ is a geodesic from $x_*$ to $x$ and $g$ is a geodesic from $x_*$ to $y$, we 
have $f(r)=g(r)$ for every $r\in [0,\varepsilon]$.
\end{corollary}

Informally, there is only one way of leaving $x_*$ along a geodesic. By the remarks of
the beginning of this section, the same property holds for a typical point chosen 
according to the volume measure on ${\bf m}_\infty$. 

\section{Infinite volume limits and the Brownian plane}
\label{secBP}

\subsection{UIPT and UIPQ}

Theorem \ref{mainresult} deals with the convergence of uniformly distributed 
rooted $p$-angula\-tions with $n$ faces in the Gromov-Hausdorff sense, provided that
the graph distance is rescaled by the factor $n^{-1/4}$ when $n$ tends to infinity. 
Note that this rescaling is necessary if we want to obtain a compact limit. On the
other hand, one may also consider the convergence of the same random
objects without rescaling, but then in a different sense than the Gromov-Hausdorff
convergence. This leads to infinite random lattices.

We consider possibly infinite (multi)graphs that are always connected, pointed (meaning that there is a distinguished
vertex $\rho$ called the root vertex) and locally finite in the sense that the degree of every vertex is
finite. As previously, these graphs are equipped with the graph distance. A ball of radius
$k$ centered at a vertex $v$ of the graph is then viewed as the subgraph consisting
of all vertices at distance less than or equal to $k$ from $v$ and the edges connecting these vertices.
A sequence $(G_n,\rho_n)$
of pointed graphs is said to converge locally to a limiting pointed graph $(G,\rho)$ if for every integer $k\geq 0$, for every
$n$ sufficiently large, the
ball of radius $k$ centered at $\rho_n$ in $G_n$ is equal to the ball of radius $k$ 
centered at $\rho$ in $G$ (equality here is in the sense
of isomorphism between finite graphs with a distinguished vertex).

The following result is due to Angel and Schramm \cite{AS} in the case of triangulations
and to Krikun \cite{Kr} (see also \cite{CMM}) in the case of quadrangulations.

\begin{theorem}
\label{infinite-lattice}
Let $p=3$ or $p=4$, and  let $M_n$ be uniformly distributed over $\mathcal{M}^p_n$.
There exists a random infinite graph $M^{(p)}_\infty$ such that
$$M_n\build{\longrightarrow}_{n\to\infty}^{\rm(d)} M^{(p)}_\infty,$$
where the convergence holds in distribution in the sense of the local convergence of graphs.
\end{theorem}

In this theorem, the convergence just means that, for every integer $k\geq 0$, the probability that the ball
of radius $k$ centered at the root vertex in $M_n$ is equal to a given graph converges to the same probability for the limit $M^{(p)}_\infty$. 
A completely different construction of $M^{(4)}_\infty$ based on a version of the CVS bijection for infinite trees was
given in \cite{CD} (see \cite{Men} for a proof of the fact that the two constructions give the same object). 

\smallskip
\noindent{\bf Remark}. One may also define the limit in Theorem \ref{infinite-lattice} as an infinite planar map,
that is, with a given embedding in the sphere, and this makes it possible to give a stronger
form of the local convergence (see \cite{AS} and \cite{Kr}). Here for simplicity we avoid dealing with
infinite planar maps, since only the graph structure will play a role in the subsequent statements.

\smallskip
The infinite random graph $M^{(p)}_\infty$ is called   
the uniform infinite planar triangulation (UIPT) when $p=3$ and the
uniform infinite planar quadrangulation (UIPQ) when $p=4$. 
Properties of these infinite random graphs have been investigated in
detail in the recent years. In particular, Angel \cite{Ang} has studied
percolation on the UIPT. The recurrence of simple random walk on the UIPT or the UIPQ 
has been obtained recently by Gurel-Gurevich and Nachmias \cite{GN}. See also 
Benjamini and Curien \cite{BC} for a proof of subdiffusivity of simple random walk on the UIPQ.

\subsection{Convergence to the Brownian plane} The following results are taken from \cite{CLG}. 
To simplify notation, we write $Q_\infty=M^{(4)}_\infty$ for  the UIPQ, which was introduced in the previous subsection. The next theorem
shows that, if we rescale the graph distance on the UIPQ by a factor tending to $0$, the resulting 
metric spaces converge (in a suitable sense) toward a limiting random non-compact metric space,
which is called the Brownian plane.  We recall that a pointed metric space is just a metric space equipped 
with a distinguished point.

\begin{theorem}
\label{conv-BPlane}
Let $V(Q_\infty)$ denote the vertex set of $Q_\infty$, which is equipped with the 
graph distance $d_{\rm gr}$, and let $\rho_{Q_\infty}$ stand for the root vertex of $Q_\infty$. There exists
a random non-compact pointed metric space $(\mathcal{P}, D_\infty,\rho_\infty)$ called the Brownian plane such that
$$(V(Q_\infty), \lambda\,d_{\rm gr}, \rho_{Q_\infty})
\build{\longrightarrow}_{\lambda\to 0}^{\rm(d)}
(\mathcal{P}, D_\infty,\rho_\infty),$$
where the convergence holds in distribution in the local Gromov-Hausdorff sense.
\end{theorem}

We refer to \cite{BBI} for a precise definition of the local Gromov-Hausdorff convergence for 
(non-compact) pointed metric spaces. In the present setting, this means that, for every
$r>0$, the ball of radius $r$ centered at $\rho_{Q_\infty}$ in $(V(Q_\infty), \lambda\,d_{\rm gr})$ will
converge in distribution in the Gromov-Hausdorff sense to the ball of radius $r$ centered at $\rho_\infty$ in the
limiting space $(\mathcal{P}, D_\infty)$ -- in fact we should use here the Gromov-Hausdorff distance
between pointed compact metric spaces, which is defined by a minor modification of
the definition in subsection \ref{GHsubsec}.

The Brownian plane $(\mathcal{P},D_\infty,\rho_\infty)$ appears in several other limit theorems. In particular,
if in the convergence of Theorem \ref{mainresult} for $p=4$ one replaces the scaling factor $n^{-1/4}$ by a
function $\beta(n)$ tending to $0$ but such that  $n^{1/4}\beta(n)$ tends to infinity, the convergence still holds in the local Gromov-Hausdorff sense
and the limit is now the Brownian plane.

Alternatively, the Brownian plane
can be viewed as the tangent cone in distribution of the Brownian map at a distinguished vertex $U$ chosen at random according to the volume measure $\Lambda$.  In fact, a stronger property holds: One can construct, on the same
probability space, both the Brownian map ${\bf m}_\infty$ and the Brownian plane $\mathcal{P}$, in such a way that, a.s., there exists
a (random) $\varepsilon >0$ such that the balls of radius $\varepsilon$
centered at the distinguished point in ${\bf m}_\infty$ and in $\mathcal{P}$ are isometric. The latter fact, together with the scale
invariance property of $\mathcal{P}$ (see subsection \ref{propBP} below), allows one
to derive many properties of the Brownian plane from those known for the Brownian map.

We will now give a precise construction of the Brownian plane. Not surprisingly, this construction
is very similar to the one developed above for the Brownian map. We consider two independent
three-dimensional Bessel processes $R$ and $R'$ started from $0$ (see e.g.~\cite{RY} for basic facts about Bessel processes). We then 
define a process $Y=(Y_t)_{t\in \mathbb{R}}$ indexed by the real line, by setting
$$Y_t:=\left\{\begin{array}{ll} R_t \qquad&\hbox{if }t\geq 0,\\
\noalign{\smallskip}
R'_{-t}\qquad&\hbox{if }t\leq 0.
\end{array}
\right.
$$
Then, for every $s,t\in\mathbb{R}$, we set
$$m_Y(s,t):= \inf_{r\in \overline{st}} Y_r,$$
with the notation $\overline{st} = [s\wedge t,s\vee t]$ if $st\geq 0$, $\overline{st} = (-\infty,s\wedge t] \cup [s\vee t,\infty)$
if $st<0$.
We define a random pseudo-distance on $\mathbb{R}$ by 
$$d_Y(s,t):= Y_s + Y_t -2\,m_Y(s,t)$$
and set $s\sim_Y t$ if $d_Y(s,t)=0$. The quotient space $\mathcal{T}_\infty=\mathbb{R}\,/\!\sim_Y$
equipped with $d_Y$ is a (non-compact) random real tree, which is
sometimes called the infinite Brownian tree.  
We write $p_\infty:\mathbb{R}\longrightarrow \mathcal{T}_\infty$
for the canonical projection and 
set $\rho_\infty:= p_\infty(0)$, which plays the role of the root of $\mathcal{T}_\infty$.
The volume measure on $\mathcal{T}_\infty$ is the image of Lebesgue measure on $\mathbb{R}$ under $p_\infty$.

We next consider Brownian motion indexed by $\mathcal{T}_\infty$. Formally, we consider
a real-valued process $(Z^\infty_t)_{t\in \mathbb{R}}$ such that,
conditionally given the process $Y$, $Z^\infty$ is a centered Gaussian process  with
covariance
$$E[Z^\infty_sZ^\infty_t\mid Y] =m_Y(s,t),$$
so that we have $Z^\infty_0=0$ and $E[(Z^\infty_s-Z^\infty_t)^2\mid Y] = d_Y(s,t)$. It is not hard to verify that
the process $Z^\infty$ has a modification with continuous paths. Then a.s. we have $Z^\infty_s=Z^\infty_t$ for every $s,t\in \mathbb{R}$
such that $d_Y(s,t)=0$ and therefore we may view $Z^\infty$ as indexed by $\mathcal{T}_\infty$.

For every $s,t\in \mathbb{R}$, we set
$$D^\circ_\infty (s,t) := Z^\infty_s + Z^\infty_t - 2\,\min_{r\in[s\wedge t,s\vee t]} Z^\infty_r.$$
We extend the definition of $D^\circ_\infty$ to $\mathcal{T}_\infty\times \mathcal{T}_\infty$
by setting for $a,b\in\mathcal{T}_\infty$,
$$D^\circ_\infty (a,b):= \min\{ D^\circ_\infty (s,t): s,t\in\mathbb{R},\; p_\infty(s)=a, p_\infty(t)=b\}.$$
Finally, we set, for every $a,b\in\mathcal{T}_\infty$,
$$D_\infty(a,b) := \inf_{a_0=a,a_1,\ldots,a_k=b} \sum_{i=1}^k D^\circ_\infty(a_{i-1},a_i)$$
where the infimum is over all choices of the integer $k\geq 1$ and of the
finite sequence $a_0,a_1,\ldots,a_k$ in $\mathcal{T}_\infty$ such that $a_0=a$ and
$a_k=b$. Then $D_\infty$ is a pseudo-distance on $\mathcal{T}_\infty$, and 
we set $a\approx b$ if $D_\infty (a,b)=0$ (this can be proved to be equivalent 
to the property $D^\circ_\infty (a,b)=0$). The Brownian plane is the quotient space $\mathcal{P}=\mathcal{T}_\infty /\approx$, which is
equipped with the metric induced by $D_\infty$ and with the distinguished 
point which is the equivalence class of $\rho_\infty$ (without risk of confusion, we still write 
$\rho_\infty$ for this equivalence class). The volume measure on $\mathcal{P}$
is the image of the volume measure on $\mathcal{T}_\infty$ under the canonical projection.

\subsection{Properties of the Brownian plane}
\label{propBP}

The Brownian plane is scale invariant, meaning that, for every $\lambda >0$, the space $(\mathcal{P},\lambda\,D_\infty,\rho_\infty)$ 
has the same distribution as $(\mathcal{P},D_\infty,\rho_\infty)$. This property can be derived from the explicit construction given above,
or from the fact that the Brownian plane is a tangent cone in distribution to the Brownian map. 

Other properties of the Brownian plane are very similar to those of the Brownian map. We state an analog 
of Proposition \ref{proper}.
\begin{proposition}
\label{proper2}
Almost surely, the Brownian plane $(\mathcal{P},D_\infty)$ is homeomorphic to the Euclidean plane $\mathbb{R}^2$ and has Hausdorff dimension $4$.
\end{proposition}
From the construction of the previous subsection, the points of the Brownian plane inherit labels $Z^\infty_x$ from the
corresponding labels on the infinite Brownian tree $\mathcal{T}_\infty$. These labels can be interpreted as  ``distances from
infinity'' in the following sense. For every $x,y\in \mathcal{P}$,
$$Z^\infty_x - Z^\infty_y=\lim_{z\to\infty} (D_\infty(x,z)-D_\infty(y,z)).$$
The existence of the preceding limit is related to a property of confluence of geodesic rays in the Brownian plane.
Recall that a geodesic ray is a continous path $\gamma:[0,\infty)\longrightarrow \mathcal{P}$ such that
$D_\infty(\gamma(t),\gamma(t'))=|t-t'|$ for every $t,t'\geq 0$.

\begin{proposition}
\label{confluplane}
Let $\gamma$ and $\gamma'$ be two geodesic rays in $\mathcal{P}$. Then there exists
two reals $\alpha$ and $\beta\geq |\alpha|$ such that $\gamma(t)=\gamma'(\alpha+t)$ for every $t\geq \beta$.
\end{proposition}

This is a kind of version at infinity of Corollary \ref{confluence}. Note that an exact analog of 
Corollary \ref{confluence} also holds for the Brownian plane, and is easy to prove by the coupling
argument mentioned above. 

The Brownian plane seems to be more tractable  than the Brownian map for explicit calculations, partly because 
of its scale invariance. Let us discuss some recent results from \cite{CLG2} that shed light on the 
probabilistic structure of the Brownian plane. For every
$r>0$, we let $B_r$ be the closed ball of radius $r$ centered at $\rho_\infty$ in $\mathcal{P}$, and we define the ``extended ball'' $B^\bullet_r$ as the complement of the unbounded connected component
of $(B_r)^c$. Informally, $B^\bullet_r$ is obtained 
by ``filling the holes'' in $B_r$. We write $| B^\bullet_r|$ for the
volume of $B^\bullet_r$.

\begin{proposition}
\label{extendeball}
For every $\lambda >0$, 
$$E[\exp(-\lambda | B^\bullet_r|)] = \frac{3^{3/2}\,\cosh((2\lambda)^{1/4}r)}{\Big(\cosh^2 ((2\lambda)^{1/4}r)+2\Big)^{3/2}}.$$
\end{proposition}

In fact, one can give a simple description of the whole process $(| B^\bullet_r|)_{r\geq 0}$.
Set $\psi(\lambda):= (\frac{8}{3})^{1/2}\lambda^{3/2}$ and recall that the continuous-state branching process with
branching mechanism $\psi$ is the Markov process in $\mathbb{R}_+$ whose transition kernels are
characterized by the Laplace transform
$$E[\exp(-\lambda X_t)\mid X_0=x]= \exp(-x\,u_t(\lambda))$$
where the function $u_t(\lambda)$ solves the differential equation $\frac{du_t(\lambda)}{dt}=-\psi(u_t(\lambda))$
with initial condition $u_0(\lambda)=\lambda$ (see e.g.~\cite{LZurich}). Note that the sample paths of $X$ are right-continuous with
left limits and that $X$ is absorbed at $0$. 
We can define a process $X^\infty=(X^\infty_t)_{t\leq 0}$, which is indexed by negative times and corresponds to $X$ ``started
from $+\infty$ at time $-\infty$ and conditioned to hit $0$ at time $0$'' (formally, we start $X$ at time $0$ from $x>0$, we then shift time
so that the hitting time of $0$ becomes $0$, and we finally let $x$ tend to $+\infty$). 

\begin{proposition}
\label{extendeballprocess}
The process $(| B^\bullet_r|)_{r\geq 0}$ has the same distribution as the process $(W_r)_{r\geq 0}$ defined by
$$W_r= \sum_{-r\leq u\leq 0} \xi_u\,(\Delta X^\infty_u)^2,$$
where $\Delta X^\infty_u$ stands for the jump of $X^\infty$ at time $u$, and, conditionally given $X^\infty$, the
nonnegative random
variables $\xi_u$ are independent and identically distributed with density
$$\frac{1}{\sqrt{2\pi}}\,x^{-5/2}\,e^{-1/2x}.$$
\end{proposition}

Informally, each jump of the process $(| B^\bullet_r|)_{r\geq 0}$ corresponds to the creation of a new connected 
component of $(B_r)^c$ (there are many such components, as shown by Proposition \ref{cactus} in the 
Brownian map case), noting that this newly created connected component will be ``swallowed'' by the
extended ball. For $r\geq 0$, the random variable $X^\infty_{-r}$ should be interpreted as the 
length, in a generalized sense, of the boundary of $B^\bullet_r$. Note that, when a connected component is
swallowed, the length of the boundary of $B^\bullet_r$ has a negative jump. The distribution of the variables $\xi_u$
then corresponds to the law of the volume of a connected component given that its boundary has length $1$.

The preceding interpretation is closely related to some results of Krikun \cite{Kr} in the discrete setting
of the UIPQ. The description of the process $(| B^\bullet_r|)_{r\geq 0}$ can also be interpreted in terms of the
``peeling process'' of the UIPT studied in \cite{Ang}.

\section{Canonical embeddings and open questions}
\label{secopen}

In this section, we come back to the questions that were discussed at the beginning of the introduction above. We would like to
have a ``canonical'' construction of a random metric $\Delta$ on the sphere $\mathbb{S}^2$, in such a way that
$$({\bf m}_\infty, D) \build{=}_{}^{\rm(d)} (\mathbb{S}^2, \Delta).$$
Furthermore we expect $\Delta$ to behave well under the conformal transformations of the sphere.

Recall that the embedding of a planar map is defined up to orientation-preser\-ving homeomorphisms of the sphere.
However, there are (almost) canonical ways of choosing these embeddings. Consider the case of simple
triangulations, that is, triangulations without loops or multiple edges. According the circle packing theorem, any such triangulation can be represented
via a circle packing of the sphere, in such a way that the vertex set of the triangulation is the set of centers of all circles, and two
vertices are linked by an edge if and only if the associated circles are tangent (see Fig.~5 for an example). This representation is in fact unique up 
to the conformal transformations of the sphere (the M\"obius transformations). 

Next suppose that, for every even integer $n\geq 2$, we have constructed a circle-packing embedding $\mathcal{C}_n$ of a 
uniformly distributed simple triangulation with $n$ faces. Write $V(\mathcal{C}_n)$ for the vertex set 
of $\mathcal{C}_n$ and $d^n_{\rm gr}$ for the graph distance on $V(\mathcal{C}_n)$. By Theorem
\ref{mainresult}, or more precisely by the extension of Theorem
\ref{mainresult} to simple triangulations found in \cite{AA}, we have
$$\Big(V(\mathcal{C}_n),(\frac{3}{2})^{1/4}n^{-1/4}\, d^n_{\rm gr}\Big)\build{\longrightarrow}_{n\to\infty}^{\rm(d)} ({\bf m}_\infty, D)$$
in the Gromov-Hausdorff sense. Note that the constant $(3/2)^{1/4}$ is different from the constant
$6^{1/4}$ in \eqref{convtri} or in Theorem
\ref{mainresult}, because we are dealing with simple triangulations. 

\begin{figure}
\begin{center}
\includegraphics[width=8cm,height=8cm]{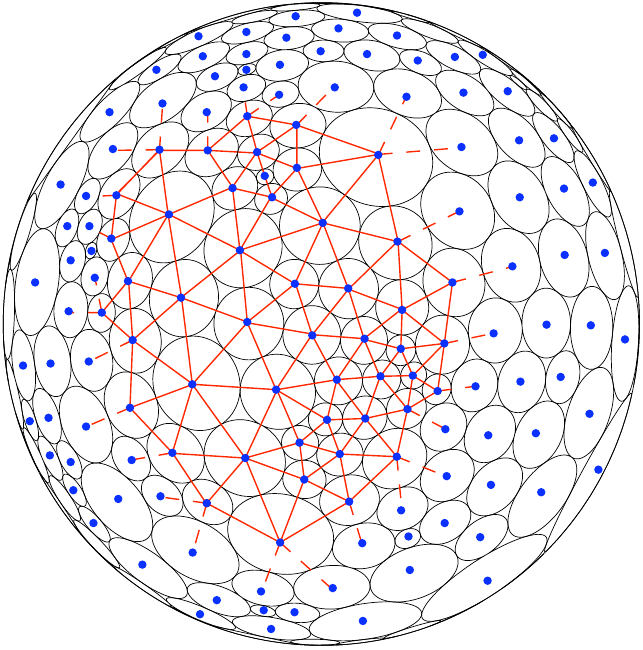}
\end{center}
\caption{A circle packing of the sphere associated with a simple triangulation.}
\end{figure}

\medskip
\noindent{\bf Conjecture.} {\it One can construct the circle packing embeddings $\mathcal{C}_n$ in such a way that 
$$\sup_{x\in\mathbb{S}^2}\Big(\min_{y\in V(\mathcal{C}_n)} |x-y|\Big) \build{\longrightarrow}_{n\to\infty}^{} 0$$
in probability, and there exists a continuous random process $(\Delta(x,y))_{x,y\in\mathbb{S}^2}$, which
is nonzero outside the diagonal and such that 
$$\sup_{x,y\in V(\mathcal{C}_n)} \Big| \Delta(x,y) - (\frac{3}{2})^{1/4} n^{-1/4}\,d^n_{\rm gr}(x,y)\Big| \build{\longrightarrow}_{n\to\infty}^{} 0$$
in probability.}

\medskip
If the conjecture holds, then it is not hard to verify that $\Delta$ defines a random metric on the sphere and that
$(\mathbb{S}^2,\Delta)$ gives a representation of the Brownian map. We note that there are other ways of defining
canonical embeddings of planar maps, which would lead to other versions of the conjecture. In particular, we can associate with
a triangulation of the sphere a Riemann surface, which is obtained by viewing each face as an equilateral triangle with sides of
length $1$ and then gluing two adjacent triangles along their common edge. The resulting Riemann surface is 
homeomorphic to the sphere and thus the uniformization theorem yields a canonical embedding, again modulo
M\"obius transformations. See \cite{Cu} for an
application of these ideas to the conformal structure of random planar maps.

As a final remark, we note that a very recent work of Miller and Sheffield \cite{MS} has introduced a new growth
process called Quantum Loewner Evolution (QLE), which might provide a direct construction of the Brownian
map, or perhaps of the Brownian plane discussed in the previous section, with conformal invariance properties.
Assuming that this construction goes through, one may expect new fascinating connections between the Brownian
map on one hand, the Schramm-Loewner evolutions and the two-dimensional Gaussian free field on the
other hand.

\end{document}